\documentclass[11pt,oneside]{amsart}
\usepackage{epsfig,psfrag,subfigure}
\usepackage{amsmath,amssymb,amsfonts,latexsym, mathrsfs}
\usepackage{graphicx,graphics,cite}
\usepackage{wrapfig}
\usepackage[usenames]{color}
\usepackage{geometry}                
\usepackage{colortbl,slashbox,booktabs,ctable,multirow}
\usepackage[table]{xcolor}
\usepackage{epstopdf}

\topmargin by -\headsep \textheight 8.9in \oddsidemargin 0.15in
\makeatletter \@addtoreset{figure}{section}
\def\thefigure{\thesection.\@arabic\c@figure} \def\fps@figure{h, t}
\@addtoreset{table}{section}
\def\thetable{\thesection.\@arabic\c@table}\def\fps@table{h, t}
\@addtoreset{equation}{section}

\makeatother

\newtheorem{algorithm}{Algorithm}[section]

\newcommand{\beq}{\begin{equation}}
\newcommand{\eeq}{\end{equation}}
\newcommand{\beqa}{\begin{eqnarray}}

\newcommand{\eeqa}{\end{eqnarray}}

\newcounter{nfig}

\begin{document}

\title{An introduction to ML($n$)BiCGStab}
\author{Man-Chung Yeung\\ \\
Department of Mathematics, University of Wyoming,
USA}
\date{}
\maketitle


\begin{abstract}
ML($n$)BiCGStab is a Krylov subspace method for the solution of large, sparse and non-symmetric
linear systems.
In theory, it is a method that lies between the well-known BiCGStab and GMRES/FOM. In fact, when $n = 1$, ML($1$)BiCGStab is BiCGStab and when $n = N$, ML($N$)BiCGStab is GMRES/FOM where $N$ is the size of the linear system. Therefore, ML($n$)BiCGStab is a bridge that connects the Lanczos-based BiCGStab and the Arnoldi-based GMRES/FOM.
In computation, ML($n$)BiCGStab can be much more stable and converge much faster than BiCGStab when
a problem with ill-condition is solved. We have tested ML($n$)BiCGStab on the standard oil reservoir simulation test data called SPE9 and found that ML($n$)BiCGStab reduced
the total computational time by more than $60\%$ when compared to BiCGStab. Tests made on the data from Matrix Market
also support the superiority of ML($n$)Bi- CGStab over BiCGStab. Because of the $O(N^2)$ storage
requirement in the full GMRES, one has to adopt a restart strategy to get the storage under control when GMRES is implemented.
In comparison, ML($n$)BiCGStab is a method with only $O(nN)$ storage requirement and therefore it does not need a restart strategy. In this paper, we
introduce ML($n$)BiCGStab (in particular, a new algorithm involving $\bf A$-transpose), its relations to some existing methods and its implementations.
\end{abstract}




\section{Introduction}
\label{Introduction}
ML($n$)BiCGStab is a Krylov subspace method for the solution of the linear system
\begin{equation} \label{equ:7-9-3}
{\bf A} {\bf x} = {\bf b}
\end{equation}
where ${\bf A} \in {\mathbb C}^{N \times N}$
and ${\bf b} \in {\mathbb C}^N$. It was introduced by Yeung and Chan\cite{yeungchan} in 1999 and its algorithms were recently reformulated by Yeung\cite{yeung}. ML($n$)BiCGStab is a natural generalization of BiCGStab by van der Vorst\cite{van}, built on the multiple starting Lanczos process rather than on the single starting Lanczos process. Its derivation relies on the techniques introduced by Sonneveld\cite{sonn} and van der Vorst\cite{van} in the construction of CGS and BiCGStab.
There have been three algorithms associated with the ML($n$)BiCGStab method so far, depending on how the residual vector ${\bf r}_k$ is defined and whether or not the Hermitian transpose ${\bf A}^H$ is used. In this paper, we shall simply introduce the algorithms and address some implementation issues. For more detailed, one is referred to \cite{yeung}.

Other
extensions
of BiCGStab exist. Among them are BiCGStab2 by Gutknecht\cite{gut}, BiCGStab($l$) by Sleijpen and Fokkema\cite{SF} and CPBi-CG by Zhang\cite{zhang}.

The outline of the paper is as follows. In \S\ref{sec:8-1}, we
introduce index functions which are helpful in presenting the ML($n$)BiCGStab algorithms. In \S\ref{sec:12-24-1}, we present the
ML($n$)BiCG algorithm from \cite{yeungchan}, from which
ML($n$)BiCGStab algorithms were derived.
In \S\ref{sec:12-24-2}, we introduce the ML($n$)BiCGStab algorithms and their relationships
with some existing methods. In \S\ref{sec:12-19}, implementation
issues are addressed and conclusions are made in \S\ref{sec:conc}.

\section{Index Functions}\label{sec:8-1} Let be given a positive integer $n$.
For all integers $k$, we define
$$
\begin{array}{lcl}
g_n(k) = \lfloor ( k - 1 ) /n \rfloor & \mbox{and}& r_n(k) = k - n
g_n(k)
\end{array}
$$
where $\lfloor \,
\cdot \,\rfloor$ rounds its argument to the nearest integer towards
minus infinity. We call $g_n$ and $r_n$ index functions; they are
defined on ${\mathbb Z}$, the set of all integers, with ranges ${\mathbb Z}$ and $\{1, 2, \cdots,
n\}$, respectively.

\begin{table}
\begin{center}%
\caption{Simple illustration of the index functions for $n = 3$.}
\begin{tabular}[t]{c|rlllcl}
 $k$ &0
  & 1 2 3 & 4 5 6 & 7 8 9 &  10 11 12 & $\cdots$ \\ \hline
  $g_n(k)$ &-1
  & 0 0 0 & 1 1 1 & 2 2 2 & 3 3 3 & $\cdots$ \\
  $r_n(k)$&3
  & 1 2 3 & 1 2 3 & 1 2 3 & 1 2 3 & $\cdots$ \\[1.0ex]
\end{tabular}%
\label{fig:7-9-1}
\end{center}%
\end{table}

If we write
\begin{equation}\label{equ:8-1}
k = j n + i
\end{equation}
with $1 \leq i \leq n$ and $j \in {\mathbb Z}$, then
$$
\begin{array}{rcl} g_n (j n + i) = j &\mbox{and}& r_n(j n
+ i) = i.
\end{array}
$$
 Table \ref{fig:7-9-1} illustrates the behavior of $g_n$ and $r_n$
with $n=3$.

\section{A ML($n$)BiCG Algorithm} \label{sec:12-24-1} Parallel to the derivation of
BiCGStab from BiCG by Fletcher\cite{fletcher}, ML($n$)BiCGStab was derived from a BiCG-like
method named ML($n$)BiCG, which was constructed based on the multiple starting Lanczos process with $n$ left starting vectors and a
single right starting vector.

Let be given $n$ vectors ${\bf q}_1, \ldots, {\bf
q}_n \in {\mathbb C}^N$, which we call {\it left starting vectors} or {\it shadow vectors}. Set
\begin{equation}
\begin{array}{lll}
 {\bf p}_{k} = \left( {\bf A}^H \right)^{g_n(k)}
{\bf q}_{r_n(k)}, & & k = 1, 2, 3, \cdots.
\end{array}
\label{equ:7-9-5}
\end{equation}
The following algorithm for the solution of
eqn (\ref{equ:7-9-3}) is from
\cite{yeungchan}.\\

\begin{algorithm}{\bf ML($n$)BiCG}
\label{alg:1} \vspace{.2cm}
\begin{tabbing}
x\=xxx\= xxx\=xxx\=xxx\=xxx\=xxx\kill \>1. \> Choose an initial
guess $\widehat{\bf x}_0$ and $n$ vectors ${\bf q}_1, {\bf q}_2,
\cdots,
{\bf q}_n$. \\
\>2. \>  Compute $\widehat{\bf r}_0 = {\bf b} - {\bf A} \widehat{\bf
x}_0$ and
set ${\bf p}_1 = {\bf q}_1$, $\widehat{\bf g}_0 = \widehat{\bf r}_0$. \\
\>3. \>For $k = 1, 2, 3, \cdots$, until convergence: \\
\>4. \>\>$\alpha_k = {\bf p}_k^H \widehat{\bf r}_{k-1} / {\bf p}_k^H
{\bf A} \widehat{\bf g}_{k-1}$; \\
\>5. \>\>$\widehat{\bf x}_k = \widehat{\bf x}_{k-1} + \alpha_k \widehat{\bf g}_{k-1}$; \\
\>6. \>\> $\widehat{\bf r}_k = \widehat{\bf r}_{k-1} - \alpha_k {\bf A} \widehat{\bf g}_{k-1}$; \\
\>7. \>\>For $s = \max (k - n, 0), \cdots, k - 1$ \\
\>8. \>\>\>$\beta^{(k)}_{s} = - {\bf p}^H_{s+1} {\bf A}
 \left(\widehat{\bf r}_k + \sum_{t = \max (k - n, 0) }^{s-1} \beta^{(k)}_t \widehat{\bf g}_t \right)
\big/{\bf p}^H_{s+1} {\bf A} \widehat{\bf g}_s$; \\
\>9. \>\>End \\
\>10.\>\> $\widehat{\bf g}_k = \widehat{\bf r}_k + \sum_{s = \max (k
- n, 0) }^{k -1} \beta_{s}^{(k)}
            \widehat{\bf g}_{s}$; \\
\>11. \>\>Compute ${\bf p}_{k+1}$ according to eqn (\ref{equ:7-9-5}) \\
\>12. \> End
\end{tabbing}
\end{algorithm}
\vspace{.2cm}

Even though
the algorithm has not been
tested,
it is believed to be numerically instable because of Line 11 in
which the shadow vectors are repeatedly multiplied by ${\bf
A}^H$, a type of operation which is highly sensitive to round-off
errors. The algorithm has been introduced only for the purpose of
developing ML($n$)BiCGStab algorithms.

Relations to some other methods:
\begin{enumerate}
\item {\it Relation with FOM by Saad and Schultz\cite{saad2}}. Consider the case where $n \geq N$. If we choose ${\bf q}_k = \widehat{\bf
r}_{k-1}$ in Algorithm \ref{alg:1} (it is possible since
$\widehat{\bf r}_{k-1}$ is computed before ${\bf q}_k$ is used in Line 11),
then Algorithm
\ref{alg:1} is
a FOM algorithm.

\item {\it Relation with GMRES by Saad and Schultz\cite{saad2}}. Consider the case where $n \geq N$. If we choose ${\bf q}_k = {\bf A} \widehat{\bf
r}_{k-1}$ in Algorithm \ref{alg:1},
then Algorithm
\ref{alg:1} is
a GMRES algorithm.

\item {\it Relation with BiCG}. When $n = 1$, Algorithm
\ref{alg:1} is a
BiCG algorithm.
\end{enumerate}

\section{ML($n$)BiCGStab Algorithms} \label{sec:12-24-2} There are three algorithms for the
ML($n$)BiCGStab method. All were derived from Algorithm \ref{alg:1}. The first two algorithms do not
involve ${\bf A}^H$ in their implementation and can be found in \cite{yeung}. The third one, however, needs ${\bf A}^H$ and is new. Therefore, we spend more space here on the the third algorithm.

\subsection{First Algorithm}\label{sec:9-12}
Let $\Omega_k(\lambda)$ be the polynomial of degree $k$ defined by
$$
\Omega_k(\lambda) = \left\{ \begin{array}{lcl} 1& & \mbox{if } k = 0\\
(1-\omega_k \lambda ) \Omega_{k-1}(\lambda)& & \mbox{if } k > 0.
\end{array} \right.
$$
If we define the ML($n$)BiCGStab residual ${\bf r}_k$ by
$$
{\bf r}_k = \left\{ \begin{array}{lcl}
\Omega_{g_n(k)+1}({\bf A}) \,\widehat{\bf r}_k, & & \mbox{if } k \geq 1,\\
\widehat{\bf r}_0, & & \mbox{if } k = 0,
\end{array}\right.
$$
then Algorithm \ref{alg:1} will lead to the first ML($n$)BiCGStab algorithm (Algorithm 4.1 in \cite{yeung}).
Computational and storage cost based on its preconditioned version (Algorithm 9.1 in \cite{yeung}) is presented in Table \ref{tab:10-28-1}.

\begin{table}[tbp]
\begin{center}
\caption{Average cost per iteration of the first ML($n$)BiCGStab algorithm
and its storage.}
\begin{tabular}{|c|c|c|c|}  \hline
Preconditioning ${\bf M}^{-1} {\bf v}$ & $\displaystyle{1 +
\frac{1}{n}}$ & ${\bf u} \pm {\bf v}, \,\, \alpha {\bf v}$ &
$\displaystyle{\max (4 - \frac{5}{n}, 0)}$ \\ \hline Matvec ${\bf A} {\bf
v}$ & $\displaystyle{1 + \frac{1}{n}}$ & Saxpy ${\bf u} + \alpha
{\bf v}$ & $\displaystyle{\max (2.5 n + 0.5 + \frac{1}{n}, 6)}$
 \\ \hline
dot product $\displaystyle{ {\bf u}^H {\bf v}}$ & $\displaystyle{n
+ 1+ \frac{2}{n}}$ &Storage  &${\bf A} + {\bf M} +$
\\
 & & & $
(4 n+4) N +
O(n)$ \\ \hline
\end{tabular}
\end{center} \label{tab:10-28-1}
\end{table}

Relations to some other methods: this first algorithm is a BiCGStab algorithm when $n = 1$.

\subsection{Second Algorithm}\label{sec:12-24-3}
If we define the ML($n$)BiCGStab residual ${\bf r}_k$ by
\begin{equation} \label{equ:9-24-1}
{\bf r}_k = \left\{ \begin{array}{lcl}
\Omega_{g_n(k+1)}({\bf A}) \,\widehat{\bf r}_k, & & \mbox{if } k \geq 1,\\
\widehat{\bf r}_0, & & \mbox{if } k = 0,
\end{array}\right.
\end{equation}
then Algorithm \ref{alg:1} will lead to the second ML($n$)BiCGStab algorithm (Algorithm 5.1 in \cite{yeung}).
Computational and storage cost based on its preconditioned version (Algorithm 9.2 in \cite{yeung}) is presented in Table \ref{tab:10-29-1}.

\begin{table}[tbp]
\begin{center}
\caption{Average cost per iteration of
the second ML($n$)BiCGStab algorithm and its storage.}
\begin{tabular}{|c|c|c|c|}  \hline
Preconditioning ${\bf M}^{-1} {\bf v}$ & $\displaystyle{1 +
\frac{1}{n}}$  & ${\bf u} \pm {\bf v}, \,\, \alpha {\bf v}$ &
$\displaystyle{1}$ \\ \hline Matvec ${\bf A} {\bf v}$ &
$\displaystyle{1 + \frac{1}{n}}$  & Saxpy ${\bf u} + \alpha {\bf
v}$ & $\displaystyle{2 n + 2 + \frac{2}{n}}$  \\ \hline dot product
$\displaystyle{ {\bf u}^H {\bf v}}$ & $\displaystyle{n + 1+
\frac{2}{n}}$ & Storage &${\bf A} + {\bf M} +$\\
& & & $
(3 n+5) N +O(n)$
\\ \hline
\end{tabular}
\end{center} \label{tab:10-29-1}
\end{table}

Relations to some other methods:
\begin{enumerate}
\item {\it Relation with FOM}.
Consider the case where $n \geq N$. If we choose ${\bf q}_k
= {\bf r}_{k-1}$, then this algorithm is a FOM algorithm.

\item {\it Relation with GMRES}.
Consider the case where $n \geq N$. If we choose ${\bf q}_k
= {\bf A} {\bf r}_{k-1}$, then this algorithm is a GMRES algorithm.

\item {\it Relation with BiCGStab}. When $n = 1$, this
algorithm is a BiCGStab algorithm.

\item {\it Relation with IDR($s$) by Sonneveld and van Gijzen\cite{sonn1, gs10}}. This algorithm is a IDR($n$) algorithm.
\end{enumerate}


\subsection{Third Algorithm} If we define the ML($n$)BiCGStab residual ${\bf r}_k$ by eqn (\ref{equ:9-24-1}) and get ${\bf A}^H$ involved in its implementation,
then through the derivation stages \#5 - \#8 in \cite{yeung}, Algorithm \ref{alg:1} will lead to the following ML($n$)BiCGStab algorithm which we name
ML($n$)BiCGStabt, standing for ML($n$)BiCGStab with ${\bf A}$-transpose.\\

\begin{algorithm}{\bf ML($n$)BiCGStabt without preconditioning 
} \label{alg:13}
\vspace{.2cm}
\begin{tabbing}
x\=xxx\= xxx\=xxx\=xxx\=xxx\=xxx\=xxx\=xxx\=xxx\kill
\>1. \> Choose an initial guess ${\bf
x}_0$ and $n$ vectors ${\bf q}_1, {\bf q}_2, \cdots,
{\bf q}_n$. \\
\>2. \> Compute $[{\bf f}_1, \cdots, {\bf f}_{n-1}] = {\bf A}^H [{\bf q}_1, \cdots, {\bf q}_{n-1}]$. \\
\>3. \>  Compute ${\bf r}_0 = {\bf b} - {\bf A} {\bf x}_0$ and ${\bf
g}_0 = {\bf r}_0,\,\, {\bf w}_0 = {\bf A}
{\bf g}_0,\,\, c_0 = {\bf q}^H_{1} {\bf w}_0$. \\
\>4. \>For $k = 1, 2, \cdots$, until convergence: \\
\>5. \>\>$\alpha_k = {\bf q}_{r_n(k)}^H {\bf r}_{k-1} / c_{k-1}$;\\
\>6.\>\> If $r_n(k) < n$\\
\>7. \>\>\> ${\bf x}_k = {\bf x}_{k-1} + \alpha_k {\bf g}_{k-1}$;
\,
${\bf r}_k = {\bf r}_{k-1} - \alpha_k {\bf w}_{k-1}$;
\\
\>8. \>\>\> ${\bf z}_w = {\bf r}_{k}; \,\, {\bf g}_k = {\bf 0}$;
\\
\>9. \>\>\>For $s = \max (k - n, 0), \cdots, g_n(k)n - 1$ \\
\>10. \>\>\>\>$\tilde{\beta}^{(k)}_{s} = - {\bf q}^H_{r_n(s+1)}
 {\bf z}_w \big/
 c_s$; \,\,\,\,\,\, \% $\tilde{\beta}^{(k)}_{s} = -\omega_{g_n(k+1)} \beta^{(k)}_{s}$\\
 \>11. \>\>\>\>${\bf z}_w = {\bf z}_w + \tilde{\beta}^{(k)}_{s} {\bf w}_s$;\\
 \>12. \>\>\>\>${\bf g}_k = {\bf g}_k + \tilde{\beta}^{(k)}_{s} {\bf g}_s$; \\
\>13. \>\>\>End
\\
\>14. \>\>\>${\bf g}_k = {\bf z}_w - \frac{1}{\omega_{g_n(k+1)}} {\bf g}_k$;
\\
\>15. \>\>\>For $s = g_n(k)n, \cdots, k - 1$ \\
\>16. \>\>\>\>$\beta^{(k)}_{s} = - {\bf f}_{r_n(s+1)}^H {\bf g}_k
\big/ c_s$; \\
\>17. \>\>\>\>${\bf g}_k = {\bf g}_k + \beta^{(k)}_{s} {\bf g}_s$; \\
\>18. \>\>\>End
\\
\>19.\>\>Else\\
\>20. \>\>\> $ {\bf x}_k = {\bf x}_{k-1} + \alpha_k {\bf g}_{k-1}$;
\\
\>21. \>\>\>
$ {\bf u}_k = {\bf r}_{k-1} - \alpha_k {\bf w}_{k-1}$;
\\
\>22. \>\>\> $\omega_{g_n(k+1)} =  ({\bf A} {\bf u}_k)^H {\bf u}_k / \|{\bf A}{\bf u}_k \|_2^2$;\\
\>23.\>\>\>${\bf x}_k = {\bf x}_k  +\omega_{g_n(k+1)} {\bf u}_k$;\,
${\bf r}_k = -\omega_{g_n(k+1)} {\bf A}{\bf u}_k +
{\bf u}_k$; \\
\>24. \>\>\> ${\bf z}_w = {\bf r}_{k}; \,\, {\bf g}_k = {\bf 0}$;
\\
\>25. \>\>\>For $s = g_n(k)n, \cdots, k - 1$ \\
\>26. \>\>\>\>$\tilde{\beta}^{(k)}_{s} = - {\bf q}^H_{r_n(s+1)} {\bf z}_w
\big/
 c_s$; \,\,\,\,\,\, \% $\tilde{\beta}^{(k)}_{s} = -\omega_{g_n(k+1)} \beta^{(k)}_{s}$\\
 \>27. \>\>\>\>${\bf z}_w = {\bf z}_w + \tilde{\beta}^{(k)}_{s} {\bf w}_s$;\\
 \>28. \>\>\>\>${\bf g}_k = {\bf g}_k + \tilde{\beta}^{(k)}_{s} {\bf g}_s$;\\
\>29. \>\>\>End \\
\>30.\>\>\> ${\bf g}_k = {\bf z}_w  - \frac{1}{\omega_{g_n(k+1)}} {\bf g}_{k}
$; \\
\>31.\>\>End\\
\>32.\>\> ${\bf w}_k = {\bf A} {\bf g}_k
$; \,
$c_k = {\bf q}_{r_n(k+1)}^H {\bf w}_k$;\\
\>33. \> End
\end{tabbing}
\end{algorithm}\vspace{.2cm}

A preconditioned version of Algorithm \ref{alg:13} can be obtained by
applying it to ${\bf A} {\bf M}^{-1} {\bf y} = {\bf b}$,
then recovering ${\bf x}$ through
${\bf x} = {\bf M}^{-1} {\bf y}$.
The resulting preconditioned algorithm and its Matlab code are attached in \S\ref{sec:apen}.
Computational and storage cost
is presented in Table \ref{tab:5-8-1}.

\begin{table}[tbp]\label{tab:5-8-1}
\begin{center}
\caption{Average cost per iteration of
preconditioned ML($n$)BiCGStabt and its storage. This table does not count the cost in Lines 1-2 of Algorithm \ref{alg:9-26-1}.}
\begin{tabular}{|c|c|c|c|}  \hline
Preconditioning ${\bf M}^{-1} {\bf v}$ & $\displaystyle{1 +
\frac{1}{n}}$  & ${\bf u} \pm {\bf v}, \,\, \alpha {\bf v}$ &
$\displaystyle{1}$ \\ \hline Matvec ${\bf A} {\bf v}$ &
$\displaystyle{1 + \frac{1}{n}}$  & Saxpy ${\bf u} + \alpha {\bf
v}$ & $\displaystyle{1.5 n + 2.5 + \frac{2}{n}}$  \\ \hline dot product
$\displaystyle{ {\bf u}^H {\bf v}}$ & $\displaystyle{n + 1+
\frac{2}{n}}$ & Storage &${\bf A} + {\bf M} +$
\\
 & & & $(4 n+4) N +O(n)$\\ \hline
\end{tabular}
\end{center}
\end{table}

Relations to some other methods: Algorithm \ref{alg:13} is a BiCGStab algorithm when $n = 1$.

%
%
%



\section{Implementation Issues}\label{sec:12-19}
%
The following test data were downloaded from
Matrix Market. 
More experiments can be found in \cite{yeung, yeungchan}.

\begin{enumerate}
\item {\it utm5940}, TOKAMAK Nuclear Physics (Plasmas). {\it utm5940}
contains a $5940 \times 5940$ real unsymmetric matrix ${\bf A}$ with
$83,842$ nonzero entries and a real right-hand side $\bf b$.

\item {\it qc2534}, H2PLUS Quantum Chemistry, NEP Collection.
{\it qc2534} contains a $2534 \times 2534$ complex symmetric
indefinite matrix with $463,360$ nonzero entries, but does not
provide the right-hand side ${\bf b}$. 
We set
${\bf b} = {\bf A} {\bf 1}$ with ${\bf 1} = [1, ,1, \cdots, 1]^T$.
\end{enumerate}

All computing in this section was done in Matlab Version 7.1 on a Windows XP
machine with a Pentium 4 processor. $ILU(0)$ preconditioners (p.294,
\cite{saad}) were used, initial guess was ${\bf x}_0 = {\bf 0}$ and the
stopping criterion was
$$
\| {\bf r}_k \|_2 / \| {\bf b}\|_2 < 10^{-7}
$$
where ${\bf r}_k$ was the computed residual. Shadow vectors ${\bf Q} = [{\bf q}_1, {\bf q}_2, \cdots,
{\bf q}_n]$ were chosen to be
${\bf Q} = [{\bf r}_0, randn(N,n-1)]$
  for {\it utm5940} and
${\bf Q} = [{\bf r}_0,
randn(N,n-1) + sqrt(-1)* randn(N,n-1)]$ for {\it qc2534}.

For the convenience of our presentation, let us introduce the
following functions:

\begin{enumerate}
\item[(a)]
$T_{conv}(n)$ is the time that a ML($n$)BiCGStab algorithm takes to
converge.



\item[(b)] $E(n)
\equiv \| {\bf b} - {\bf
A} {\bf x}\|_2 / \|{\bf b}\|_2$ is the true relative error of ${\bf
x}$ where ${\bf x}$ is the computed solution output by a
ML($n$)BiCGStab algorithm when it converges.
\end{enumerate}

\subsection{Stability}
\begin{figure}[htbp]
\centerline{\hbox{
\psfig{figure=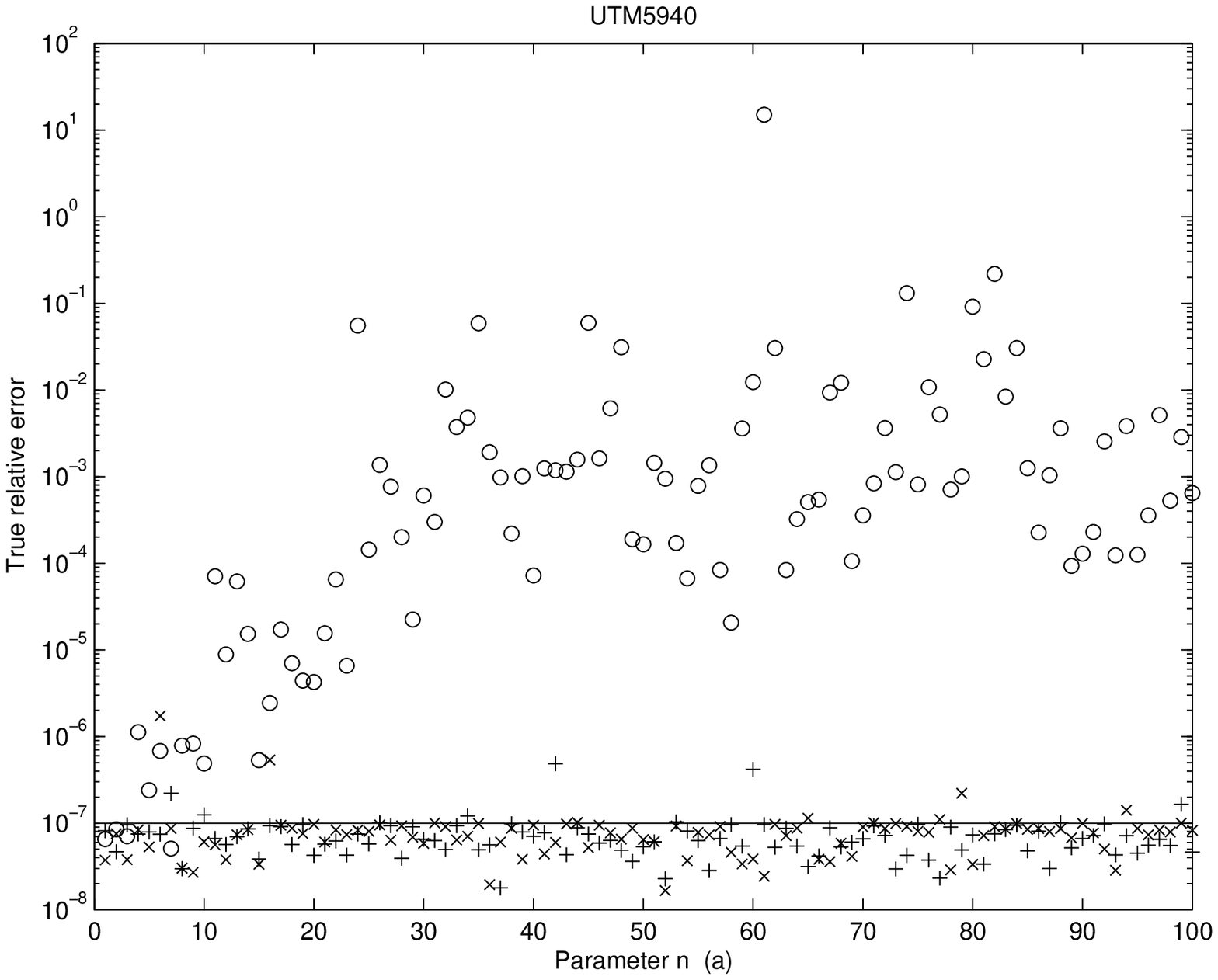,height=6.5cm}
\psfig{figure=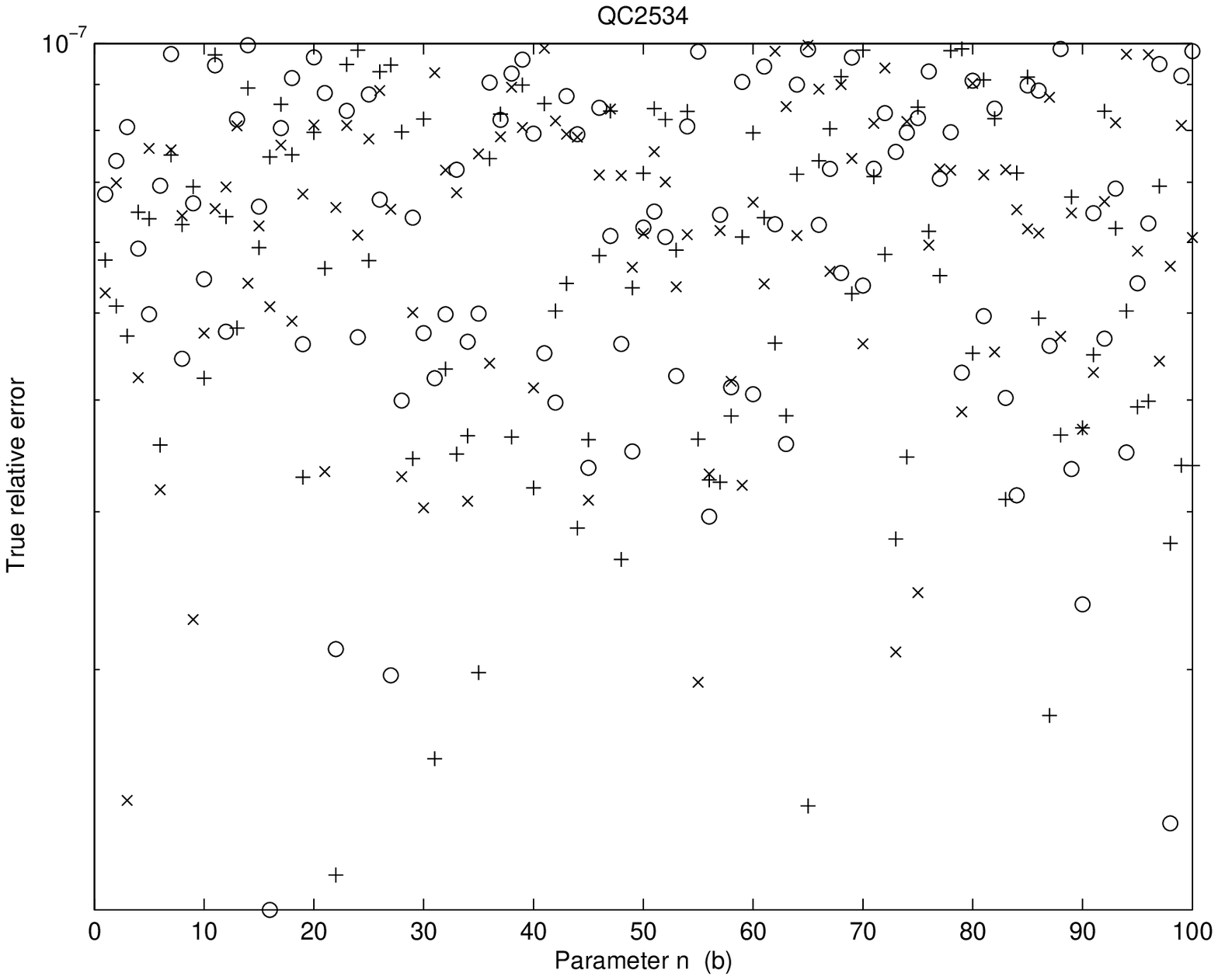,height=6.5cm}
 }} \caption{Graphs of $E(n)$ against $n$.
 First algorithm: $\times$-mark; Second algorithm: o-mark; Third algorithm: $+$-mark;
$10^{-7}$: solid line.
} \label{fig:11-24-2}
\end{figure}

The graphs of $E(n)$ are plotted in Figure \ref{fig:11-24-2}. It can be seen that
the computed ${\bf r}_k$ by the second algorithm can easily
diverges from its exact counterpart ${\bf b} - {\bf A} {\bf x}_k$.
This divergence
becomes significant when $n \geq
4$ for {\it utm5940}.
By contrast, the computed relative errors $\|{\bf r}_k\|_2 / \| {\bf
b}\|_2$ by the first and the third algorithms well approximate their
corresponding true ones.
Thus, from this point of view, we consider that the first and the third algorithms are numerically more stable than the second algorithm.

\subsection{Choice of $n$}\label{subsec:5-20}
From the experiments in \cite{yeung, yeungchan}, we have
observed that ML($n$)BiCGStab behaves more and more robust as $n$ is
increased. So, for an ill-conditioned problem, we would tend to
suggest a large $n$ for ML($n$)BiCGStab. On the other hand,
ML($n$)BiCGStab minimizes $\|{\bf r}_k\|_2$ once every $n$
iterations. The convergence of a well-conditioned problem is usually
accelerated by the minimization steps. So, when a problem is
well-conditioned, we would suggest a small $n$.

In \cite{sonn1, gs10}, it was suggested to fix $s = 4$ or $8$ for the general use of IDR($s$). This good idea also applies to ML($n$)BiCGStab, namely, fixing $n = 4$ or $8$ in its general use.

We believe that the most powerfulness of ML($n$)BiCGStab is in the solution of a sequence of linear systems.
We once tested the first algorithm (see Algorithm 9.1 in \cite{yeung}) with $n = 9$ and $\kappa = 0$ (see \S\ref{sec:6-5} for $\kappa$)
on the standard oil reservoir simulation test data called SPE9
and found that ML($n$)BiCGStab reduced the total computational time by over $70\%$ when compared to BiCGStab. A later test on SPE9 with Code \#4 in \cite{yeung} showed that a $60\%$ reduction in time can be reached.

Code \#4 is a design of automatic selection of the parameter $n$ during the solution of a sequence of linear systems. Let $t1$ and $t2$ denote the times to solve the previous and the current systems respectively. Then the
 basic idea behind Code \#4 is: if $t1 > t2$, then increase $n$ to $n + step$ when solving the next system; otherwise, decrease $n$ to $n - step$. Here $step$ is the search step size.

\begin{figure}[htbp]
\centerline{\hbox{
\psfig{figure=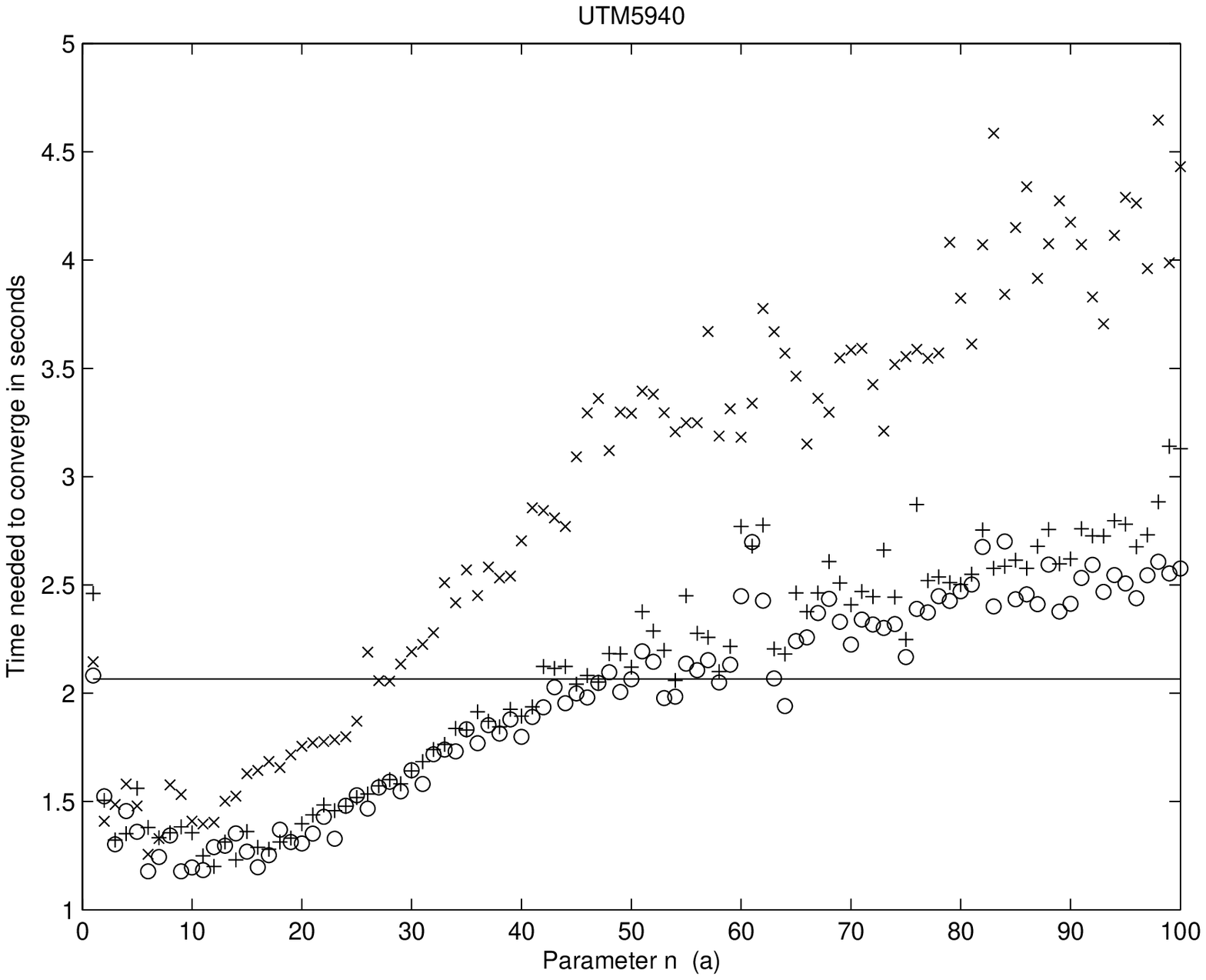,height=6.5cm}
\psfig{figure=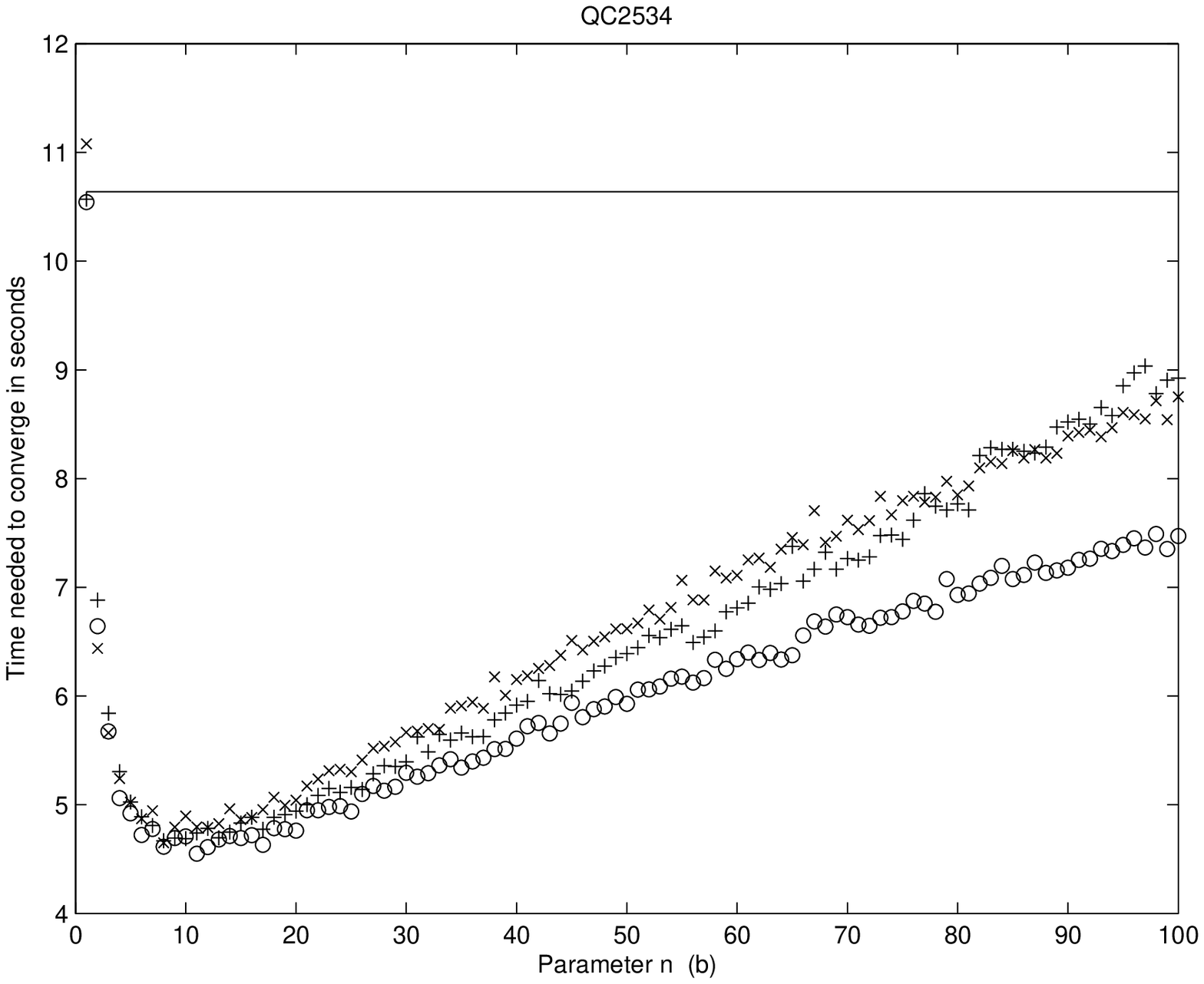,height=6.5cm}
 }} \caption{Graphs of $T_{conv}(n)$ against $n$. First algorithm: $\times$-mark; Second algorithm: o-mark; Third algorithm: $+$-mark; BiCGStab: Solid line. BiCGStab took $2.07$ and $10.64$ seconds to converge for $utm5940$ and $qc2534$ respectively.
} \label{fig:12-3-1}
\end{figure}

We also plot the graphs of $T_{conv}(n)$ in Figure \ref{fig:12-3-1}
to provide more
information on how $n$ affects the performance of ML($n$)BiCGStab.

\subsection{Choice of $\omega$}\label{sec:6-5}
The standard choice for the
$\omega_{g_n(k+1)}$ in Algorithm \ref{alg:13} (see Line 22) is
$
\omega_{g_n(k+1)} = ({\bf A} {\bf u}_{k})^H {\bf u}_{k} /
\| {\bf A} {\bf u}_{k} \|_2^2$. 
This choice of $\omega_{g_n(k+1)}$ minimizes the $2$-norm of ${\bf
r}_{k} = -\omega_{g_n(k+1)} {\bf A} {\bf u}_{k} + {\bf
u}_{k}$ (Line 23), but sometimes can cause instability due to that it can
be very small during an execution. The following  remedy to guard $\omega_{g_n(k+1)}$ away from zero has been
proposed in \cite{Svan1}:
\begin{equation}\label{equ:12-6-2}
\begin{array}{l}
 \omega_{g_n(k+1)}
= ({\bf A} {\bf u}_{k})^H {\bf u}_{k}
/ \| {\bf A} {\bf u}_{k} \|_2^2;\\
\rho = ({\bf A} {\bf u}_{k})^H {\bf u}_{k} / (\|
{\bf A} {\bf u}_{k} \|_2\,\, \| {\bf u}_{k}
\|_2);\\
\mbox{if  } |\rho| < \kappa,\,\,\, \omega_{g_n(k+1)} = \kappa \omega_{g_n(k+1)}
/ |\rho|; \,\,\, \mbox{end}
\end{array}
\end{equation}
where $\kappa$ is a user-defined parameter.
See the numerical
experiments in \cite{sonn1, yeung} for more information about eqns (\ref{equ:12-6-2}).


\section{Conclusions}\label{sec:conc}

ML($n$)BiCGStab is a powerful Krylov subspace method, especially in the solution of a sequence of linear systems with the parameter $n$ dynamically chosen (see \cite{yeung} for detail). This method has three algorithms. The first two can be found in \cite{yeung} and the third is new and is presented here as Algorithm \ref{alg:13}. The third algorithm involves ${\bf A}^H$ in its implementation and behaves as stable as the first algorithm, but converges faster than the first algorithm. Compared to the second algorithm, this third algorithm is more stable, but takes more time to converge.

\section{Appendix} \label{sec:apen}
Algorithm \ref{alg:9-26-1} below is the preconditioned version of Algorithm \ref{alg:13}. 
To avoid calling the index
functions $r_n(k)$ and $g_n(k)$ every $k$-iteration, we have split
the $k$-loop into a $i$-loop and a $j$-loop where $i, j, k$ are
related by (\ref{equ:8-1}) with $1 \leq i \leq n, 0 \leq j$.

\begin{algorithm}\label{alg:9-26-1}{\bf ML($n$)BiCGStabt with preconditioning
}
\vspace{.2cm}
\begin{tabbing}
x\=xxx\= xxx\=xxx\=xxx\=xxx\=xxx\=xxx\=xxx\=xxx\kill
\>1. \> Choose an initial guess ${\bf
x}_0$ and $n$ vectors ${\bf q}_1, {\bf q}_2, \cdots,
{\bf q}_n$. \\
\>2. \> Compute $[{\bf f}_1, \cdots, {\bf f}_{n-1}] = {\bf M}^{-H} {\bf A}^H [{\bf q}_1, \cdots, {\bf q}_{n-1}]$, ${\bf r}_0 = {\bf b} - {\bf A} {\bf x}_0$ and ${\bf
g}_0 = {\bf r}_0$.\\
\>\> Compute $\hat{\bf g}_0 = {\bf M}^{-1}
{\bf r}_0, \,\,
 {\bf w}_0 = {\bf A}
\hat{\bf g}_0,\,\, c_0 = {\bf q}^H_{1} {\bf w}_0, \,\, e_0 = {\bf q}_1^H {\bf r}_0$. \\
\>3. \>For $j = 0, 1, 2, \cdots$ \\
\>4. \>\>For $i = 1, 2, \cdots, n-1$ \\
\>5. \>\>\>$\alpha_{jn+i} = e_{jn+i-1} / c_{jn+i-1}$;\\
\>6. \>\>\> ${\bf x}_{jn+i} = {\bf x}_{jn+i-1} + \alpha_{jn+i} \hat{\bf g}_{jn+i-1}$;
\\
\>7. \>\>\> ${\bf r}_{jn+i} = {\bf r}_{jn+i-1} - \alpha_{jn+i} {\bf w}_{jn+i-1}$;
\\
\>8. \>\>\> $e_{jn+i} = {\bf q}_{i+1}^H {\bf r}_{jn+i}$;\\
\>9. \>\>\> If $j \geq 1$\\
\>10. \>\>\>\>$\tilde{\beta}^{(jn+i)}_{(j-1)n+i} = - e_{jn+i} \big/
 c_{(j-1)n+i}$; \,\,\,\,\,\, \% $\tilde{\beta}^{(jn+i)}_{(j-1)n+i} = -\omega_{j} \beta^{(jn+i)}_{(j-1)n+i}$\\
 \>11. \>\>\>\>${\bf z}_w = {\bf r}_{jn+i} + \tilde{\beta}^{(jn+i)}_{(j-1)n+i} {\bf w}_{(j-1)n+i}$;\\
 \>12. \>\>\>\>${\bf g}_{jn+i} = \tilde{\beta}^{(jn+i)}_{(j-1)n+i} {\bf g}_{(j-1)n+i}$; \\
\>13. \>\>\>\>For $s = i+1, \cdots, n- 1$ \\
\>14. \>\>\>\>\>$\tilde{\beta}^{(jn+i)}_{(j-1)n+s} = - {\bf q}^H_{s+1}
 {\bf z}_w \big/
 c_{(j-1)n+s}$; \,\,\,\,\,\, \% $\tilde{\beta}^{(jn+i)}_{(j-1)n+s} = -\omega_{j} \beta^{(jn+i)}_{(j-1)n+s}$\\
 \>15. \>\>\>\>\>${\bf z}_w = {\bf z}_w + \tilde{\beta}^{(jn+i)}_{(j-1)n+s} {\bf w}_{(j-1)n+s}$;\\
 \>16. \>\>\>\>\>${\bf g}_{jn+i} = {\bf g}_{jn+i} + \tilde{\beta}^{(jn+i)}_{(j-1)n+s} {\bf g}_{(j-1)n+s}$; \\
\>17. \>\>\>\>End
\\
\>18. \>\>\>\>$\displaystyle{{\bf g}_{jn+i} = {\bf z}_w - \frac{1}{\omega_{j}} {\bf g}_{jn+i}}$;
\\
\>19. \>\>\>\>For $s = 0, \cdots, i - 1$ \\
\>20. \>\>\>\>\>$\beta^{(jn+i)}_{jn+s} = - {\bf f}_{s+1}^H {\bf g}_{jn+i}
\big/ c_{jn+s}$; \\
\>21. \>\>\>\>\>${\bf g}_{jn+i} = {\bf g}_{jn+i} + \beta^{(jn+i)}_{jn+s} {\bf g}_{jn+s}$; \\
\>22. \>\>\>\>End
\\
\>23. \>\>\>Else
\\
\>24. \>\>\>\>$\beta^{(jn+i)}_{jn} = - {\bf f}_{1}^H {\bf r}_{jn+i}
\big/ c_{jn}$; \\
\>25. \>\>\>\>${\bf g}_{jn+i} = {\bf r}_{jn+i} + \beta^{(jn+i)}_{jn} {\bf g}_{jn}$; \\
\>26. \>\>\>\>For $s = 1, \cdots, i - 1$ \\
\>27. \>\>\>\>\>$\beta^{(jn+i)}_{jn+s} = - {\bf f}_{s+1}^H {\bf g}_{jn+i}
\big/ c_{jn+s}$; \\
\>28. \>\>\>\>\>${\bf g}_{jn+i} = {\bf g}_{jn+i} + \beta^{(jn+i)}_{jn+s} {\bf g}_{jn+s}$; \\
\>29. \>\>\>\>End
\\
\>30. \>\>\>End
\\
\>31.\>\>\> $\hat{\bf g}_{jn+i} = {\bf M}^{-1} {\bf g}_{jn+i}
$; 
${\bf w}_{jn+i} = {\bf A} \hat{\bf g}_{jn+i}
$;\\
\>32.\>\>\> $c_{jn+i} = {\bf q}_{i+1}^H {\bf w}_{jn+i}$;\\
\>33. \>\> End\\
\>34. \>\>$\alpha_{jn+n} = e_{jn+n-1} / c_{jn+n-1}$;\\
\>35. \>\> ${\bf x}_{jn+n} = {\bf x}_{jn+n-1} + \alpha_{jn+n} \hat{\bf g}_{jn+n-1}$;
\\
\>36. \>\> $ {\bf u}_{jn+n} = {\bf r}_{jn+n-1} - \alpha_{jn+n} {\bf w}_{jn+n-1}$;
\\
\>37. \>\> $ \hat{\bf u}_{jn+n} = {\bf M}^{-1} {\bf u}_{jn+n}$;
\\
\>38. \>\> $\omega_{j+1} =  ({\bf A} \hat{\bf u}_{jn+n})^H {\bf u}_{jn+n} / \|{\bf A} \hat{\bf u}_{jn+n} \|_2^2$; \\
\>39.\>\>${\bf x}_{jn+n} = {\bf x}_{jn+n}  +\omega_{j+1} \hat{\bf u}_{jn+n}$; \\
\>40. \>\>${\bf r}_{jn+n} = -\omega_{j+1} {\bf A} \hat{\bf u}_{jn+n} +
{\bf u}_{jn+n}$; \\
\>41. \>\>$e_{jn+n} = {\bf q}_1^H {\bf r}_{jn+n}$;\\
\>42. \>\>$\tilde{\beta}^{(jn+n)}_{(j-1)n+n} = - e_{jn+n}
\big/
 c_{(j-1)n+n}$; \,\,\,\,\,\, \% $\tilde{\beta}^{(jn+n)}_{(j-1)n+n} = -\omega_{j+1} \beta^{(jn+n)}_{(j-1)n+n}$\\
 \>43. \>\>${\bf z}_w = {\bf r}_{jn+n} + \tilde{\beta}^{(jn+n)}_{(j-1)n+n} {\bf w}_{(j-1)n+n}$;\\
 \>44. \>\>${\bf g}_{jn+n} = \tilde{\beta}^{(jn+n)}_{(j-1)n+n} {\bf g}_{(j-1)n+n}$;\\
\>45. \>\>For $s = 1, \cdots, n - 1$ \\
\>46. \>\>\>$\tilde{\beta}^{(jn+n)}_{jn+s} = - {\bf q}^H_{s+1} {\bf z}_w
\big/
 c_{jn+s}$; \,\,\,\,\,\, \% $\tilde{\beta}^{(jn+n)}_{s+jn} = -\omega_{j+1} \beta^{(jn+n)}_{s+jn}$\\
 \>47. \>\>\>${\bf z}_w = {\bf z}_w + \tilde{\beta}^{(jn+n)}_{jn+s} {\bf w}_{jn+s}$;\\
 \>48. \>\>\>${\bf g}_{jn+n} = {\bf g}_{jn+n} + \tilde{\beta}^{(jn+n)}_{jn+s} {\bf g}_{jn+s}$;\\
\>49. \>\>End \\
\>50.\>\> $\displaystyle{{\bf g}_{jn+n} = {\bf z}_w  - \frac{1}{\omega_{j+1}} {\bf g}_{jn+n}}
$; 
$\hat{\bf g}_{jn+n} = {\bf M}^{-1} {\bf g}_{jn+n}
$; \\
\>51.\>\> ${\bf w}_{jn+n} = {\bf A} \hat{\bf g}_{jn+n}
$; 
$c_{jn+n} = {\bf q}_{1}^H {\bf w}_{jn+n}$;\\
\>52. \> End
\end{tabbing}
\end{algorithm}

\vspace{.2cm}

{\bf Matlab code of Algorithm \ref{alg:9-26-1}}

\begin{tabbing}
x\=xxx\=
xxx\=xxx\=xxx\=xxx\=xxx\=xxx\=xxx\=xxx\=xxx\=xxx\=xxx\=xxx\kill \>1.
\>function $[x,err,iter,flag] = mlbicgstabt(A,x,b,Q,M,max\_it,tol, kappa)$\\
\>2.\\
\>3.\>\% input:\>\>\>$A$:\,\,\, N-by-N matrix. $M$: \,\,\,N-by-N preconditioner matrix. \\
\>4.\>\% \>\>\>$Q$: \, N-by-n shadow matrix $[{\bf
q}_1,\cdots,{\bf q}_n]$. $x$: initial guess.\\
\>5.\>\% \>\>\> $b$:\,\,\, right hand side vector. $max\_it$:\,\,\, maximum number
of iterations.\\
\>6.\>\%\>\>\>$tol$:\,\,\, error tolerance.\\
\>7.\>\%\>\>\>$kappa$:\>\>\> (real number) minimization step controller:\\
\>8.\>\%\>\>\>\>\>\> $kappa = 0$, standard minimization\\
\>9.\>\%\>\>\>\>\>\> $kappa > 0$, Sleijpen-van der Vorst minimization \\
\>10.\>\% output:\>\>\>$x$: solution computed. $err$: error norm. $iter$: number of iterations
performed.\\
\>11.\>\% \>\>\>$flag$:\>\>\> $= 0$, solution found to
tolerance\\
\>12.\>\% \>\>\>\>\>\>  $= 1$, no convergence given $max\_it$ iterations\\
\>13.\>\% \>\>\>\>\>\>$= -1$, breakdown. \\
\>14.\>\% storage:\>\>\> $F$: $N \times (n-1)$ matrix. $G, Q, W$: $N \times n$ matrices. $A, M$: $N \times N$ matrices.\\
\>15.\>\% \>\>\>$x, r, g\_h, z, b$: $N \times 1$ matrices. $c$: $1 \times n$ matrix.\\
\>16. \\
\>17.\>\>$N = size(A,2); \,\,n = size(Q,2)$;\\
\>18.\>\>$G = zeros(N,n);\,\, W = zeros(N,n)$;\,\,\,\,\, \%
initialize
work spaces\\
\>19.\>\>if $n > 1$, $F = zeros(N,n-1)$; end\\
\>20.\>\>$c = zeros(1,n)$;\,\,\,\,\,\,\,\,\,\,\,\,\, \% end initialization\\
\>21.\>\>\\
\>22.\>\>$iter = 0;\,\, flag = 1;\,\, bnrm2 = norm(b)$;\\
\>23.\>\>if $bnrm2 == 0.0$,\, $bnrm2 = 1.0$;\, end\\
\>24.\>\>$r = b - A*x; \,\, err = norm( r ) / bnrm2$;\\
\>25.\>\>if $err < tol$,\, $flag = 0$;\,\, return,\, end\\
\>26.\>\> \\
\>27.\>\> if $n > 1$, $F = M' \backslash (A' * Q(:,1: n-1))$; end\\
\>28.\>\>$G(:,1) = r;\,\, g\_h = M \backslash r; \,\, W(:,1) = A*g\_h; \,\, c(1) = Q(:,1)'*W(:,1)$;\\
\>29.\>\>if $c(1) == 0$,\, $flag = -1$;\,\, return,\, end \\
\>30.\>\>$e = Q(:,1)'*r$; \\
\>31.\>\>\\
\>32.\>\>  for $j = 0:max\_it$\\
\>33.\>\>\>for $i = 1:n-1$\\
\>34.\>\>\>\> $alpha = e / c(i)$;
\, $x = x + alpha*g\_h$; \,
$r = r - alpha*W(:, i)$;\\
\>35.\>\>\>\> $err = norm(r)/bnrm2$;\,
$iter = iter + 1$;\\
\>36.\>\>\>\>if $err < tol$,\, $flag = 0$;\,\, return,\, end\\
\>37.\>\>\>\>if $iter >= max\_it$,\, return,\, end\\
\>38.\>\>\>\>\\
\>39.\>\>\>\>$e = Q(:,i+1)'*r$; \\
\>40.\>\>\>\>if $j >= 1$\\
\>41.\>\>\>\>\>$beta = -e / c(i+1)$; \\
\>42.\>\>\>\>\>$W(:,i+1) = r + beta*W(:,i+1)$; \\
\>43.\>\>\>\>\>$G(:,i+1) = beta*G(:,i+1)$; \\
\>44.\>\>\>\>\> for $s = i+1 : n-1$\\
\>45.\>\>\>\>\>\>$beta = -Q(:,s+1)'*W(:,i+1)/c(s+1)$; \\
\>46.\>\>\>\>\>\>$W(:,i+1) = W(:,i+1) + beta*W(:,s+1)$;\\
\>47.\>\>\>\>\>\>$G(:,i+1) = G(:,i+1) + beta*G(:,s+1)$;\\
\>48.\>\>\>\>\>end \\
\>49.\>\>\>\>\>$G(:,i+1) = W(:,i+1) - G(:,i+1)./omega$; \\
\>50.\>\>\>\>\>for $s = 0:i-1$\\
\>51.\>\>\>\>\>\>$beta = -F(:,s+1)'*G(:,i+1) / c(s+1)$;\\
\>52.\>\>\>\>\>\>$G(:,i+1) = G(:,i+1) + beta*G(:,s+1)$; \\
\>53.\>\>\>\>\>end \\
\>54.\>\>\>\>else \\
\>55.\>\>\>\>\>$beta = -F(:,1)'*r / c(1)$; \,
$G(:,i+1) = r + beta*G(:,1)$;\\
\>56.\>\>\>\>\>for $s = 1:i-1$\\
\>57.\>\>\>\>\>\>$beta = -F(:,s+1)'*G(:,i+1) / c(s+1)$;\\
\>58.\>\>\>\>\>\>$G(:,i+1) = G(:,i+1) + beta*G(:,s+1)$; \\
\>59.\>\>\>\>\>end \\
\>60.\>\>\>\>end \\
\>61.\>\>\>\>$g\_h = M \backslash G(:,i+1)$; \,
$W(:,i+1) = A*g\_h$; \\
\>62.\>\>\>\>$c(i+1) = Q(:,i+1)'*W(:,i+1)$;\\
\>63.\>\>\>\>if $c(i+1) == 0$,\, $flag = -1$;\,\, return,\, end \\
\>64.\>\>\>end \\
\>65.\>\>\> $alpha = e / c(n)$;\,
$x = x + alpha*g\_h$; \,
$r = r - alpha*W(:,n)$; \\
\>66.\>\>\>
$err = norm(r)/bnrm2$;\\
\>67.\>\>\>if $err < tol$,\, $flag = 0; \,\, iter = iter + 1$;\, return,\, end \\
\>68.\>\>\>$g\_h = M \backslash r; \,\, z = A*g\_h$;\, $omega = z'*z$; \\
\>69.\>\>\>if $omega == 0$,\, $flag = -1$;\, return,\, end\\
\>70.\>\>\>$rho = z'*r$;\,
$omega =  rho / omega$; \\
\>71.\>\>\>if $kappa > 0$\\
\>72.\>\>\>\>$rho = rho / (norm(z)*norm(r))$;\,
$abs\_om = abs(rho)$;\\
\>73.\>\>\>\>if ($abs\_om < kappa$) \& ($abs\_om \sim = 0$)\\
\>74.\>\>\>\>\>$omega = omega*kappa/abs\_om$;\\
\>75.\>\>\>\>end\\
\>76.\>\>\>end\\
\>77.\>\>\>if $omega == 0$,\, $flag = -1$;\, return,\, end\\
\>78.\>\>\>$x = x + omega*g\_h$;\,
$r = r - omega*z$;\\
\>79.\>\>\>$err = norm(r)/bnrm2$;\,
$iter = iter + 1$;\\
\>80.\>\>\>if $err < tol$,\, $ flag = 0$;\,\, return,\, end\\
\>81.\>\>\>if $iter >= max\_it$,\, return,\, end \\
\>82.\>\>\> \\
\>83.\>\>\>$e = Q(:,1)'*r; \,\, beta = - e/c(1)$; \\
\>84.\>\>\>$W(:,1) = r + beta*W(:,1)$;\,
$G(:,1) = beta*G(:,1)$; \\
\>85.\>\>\> for $s = 1:n-1$ \\
\>86.\>\>\>\>$beta = -Q(:,s+1)'*W(:,1)/c(s+1)$;\\
\>87.\>\>\>\>$W(:,1) = W(:,1) + beta*W(:,s+1)$; \\
\>88.\>\>\>\>$G(:,1) = G(:,1) + beta*G(:,s+1)$;\\
\>89.\>\>\> end \\
\>90.\>\>\>$G(:,1) = W(:,1) - G(:,1)./omega$;\,
$g\_h = M \backslash G(:,1)$; \\
\>91.\>\>\>$W(:,1) = A * g\_h$; \,
$c(1) = Q(:,1)'*W(:,1)$;\\
\>92.\>\>\>if $c(1) == 0$,\, $flag = -1$;\,\, return,\, end \\
\>93.\>\>  end
\end{tabbing}


{\bf
A sample execution of ML($n$)BiCGstabt}
\begin{tabbing}
x\=xxx\=
xxx\=xxx\=xxx\=xxx\=xxx\=xxx\=xxx\=xxx\=xxx\=xxx\=xxx\=xxx\kill \>1.
\>$N = 100$; $A = randn(N)$; $M = randn(N)$; $b = randn(N, 1)$;\\
\>2.\>$n = 10$; $kappa = 0.7$; $tol = 10^{-7}$; $max\_it = 3*N$;\\
\>3.\>$Q = sign(randn(N, n))$; $x = zeros(N, 1)$; $Q(:,1) = b-A*x$;\\
\>4.\>$[x,err,iter,flag] = mlbicgstabt(A,x,b,Q,M,max\_it,tol, kappa)$;
\end{tabbing}


\begin{thebibliography}{10}

\bibitem{fletcher}
Fletcher, R.,
 {\it Conjugate gradient methods for indefinite
systems,} volume 506 of Lecture Notes Math., pages 73-89.
Springer-Verlag, Berlin-Heidelberg-New York, 1976.

\bibitem{saad}
Saad, Y.,
{\it Iterative methods for sparse linear systems},
2nd edition, SIAM, Philadelphia, PA, 2003.

\bibitem{saad2}
Saad, Y. \& Schultz, M.H., {\it GMRES: A generalized minimal
residual algorithm for solving nonsymmetric linear systems,} SIAM J.
Sci. Statist. Comput., 7 (1986), pp.\ 856--869.

\bibitem{SF}
Sleijpen, G.L.G. \& Fokkema, D.R., {\it BiCGSTAB($l$) for
linear equations involving unsymmetric matrices with complex
spectrum,} ETNA, 1:11-32, 1993.


\bibitem{Svan1}
Sleijpen, G.L.G. \& van der Vorst, H.A., {\it Maintaining
convergence properties of BiCGSTAB methods in finite precision
arithmetic,} Numer. Algorithms, 10 (1995), pp.\ 203--223.

\bibitem{sonn}
Sonneveld, P., {\it CGS, a fast Lanczos-type solver for
nonsymmetric linear systems,} SIAM J. Sci. Statist. Comput., 10
(1989), pp.\ 36--52.


\bibitem{sonn1}
Sonneveld P. \& van Gijzen, M., {\it IDR(s): a family of
simple and fast algorithms for solving large nonsymmetric linear
systems,} SIAM J. Sci. Comput. Vol. 31, No. 2, pp. 1035-1062.


\bibitem{van}
van der Vorst, H.A., {\it Bi-CGSTAB: A fast and smoothly
converging variant of Bi-CG for the solution of nonsymmetric linear
systems,} SIAM J. Sci. Statist. Comput., 12 (1992), pp.\ 631--644.

\bibitem{gut}
Gutknecht, M.H., {\it Variants of BiCGStab for matrices with complex
spectrum}, SIAM J. Sci. Comput. 14, 1020-1033, 1993.


\bibitem{gs10}
van Gijzen, M. \& Sonneveld, P., {\it An elegant IDR($s$)
variant that efficiently exploits bi-orthogonality properties},
Delft University of Technology, Reports of the Department of Applied
Mathematical Analysis, Report 08-21.


\bibitem{yeung}
Yeung, M., {\it ML($n$)BiCGStab: Refomulation, Analysis and Implementation}, submitted to Numerical Mathematics: Theory, Methods and Applications.
Available at http://www.uwyo.edu/mathmyeung/r17.pdf or http://arxiv.org/abs/1011.5314v1.


\bibitem{yeungchan}
Yeung, M. \& Chan, T., {\it ML($k$)BiCGSTAB: A BiCGSTAB variant
based on multiple Lanczos starting vectors}, SIAM J. Sci. Comput.,
Vol. 21, No. 4, pp.~1263-1290, 1999.

\bibitem{zhang}
Zhang, S.L., {\it GPBi-CG: Generalized product-type
methods based on Bi-CG for solving nonsymmetric linear systems},
SIAM J. Sci. Comput., 18:537-551, 1997.

\end{thebibliography}
\end{document}